\documentclass{amsart}
\usepackage[cp1251]{inputenc}
\usepackage[english,russian]{babel}
\usepackage{amsmath}
\usepackage{amssymb}
\usepackage{amsfonts}
\def\udcs{517.9} 
\newtheorem{lemma}{Лемма}

\newtheorem*{Theorem A}{Теорема A}
\newtheorem*{Theorem 1}{Теорема 1}
\newtheorem*{Theorem 2}{Теорема 2}

\def\Re{\operatorname{\mathrm Re}}

\def\Fin{\operatorname{\mathrm Fin}}

\def\Ker{\operatorname{\mathrm Ker}}

\def\art{\operatorname{\mathrm arctg}}
\def\int{\operatorname{\mathrm int}}

\begin{document}
УДК \udcs
\thispagestyle{empty}

\title[ Интерполяции рядами экспонент
\dots]{Интерполяция рядами экспонент в $H(D),$ с вещественными узлами}

\author{С.Г. Мерзляков, С.В. Попенов}

\address{Сергей Георгиевич Мерзляков,
\newline\hphantom{iii} Институт математики c ВЦ УНЦ РАН,
\newline\hphantom{iii} ул. Чернышевского, 112, 
\newline\hphantom{iii} 450008, г. Уфа, Россия}
\email{msg2000@mail.ru}

\address{Сергей Викторович Попёнов, 
\newline\hphantom{iii} Институт математики c ВЦ УНЦ РАН,
\newline\hphantom{iii} ул. Чернышевского, 112, 
\newline\hphantom{iii} 450008, г. Уфа, Россия}
\email{spopenov@gmail.com}

\thanks{\sc S.G. Merzlyakov, S.V. Popenov, 
Interpolation by means of series of exponentials in $H(D)$ with real nodes}
\thanks{\copyright \ Мерзляков С.Г., Попёнов С.В. 2014}
\thanks{\rm Работа поддержана РФФИ (грант \No 11-01-00572-а)}
\thanks{\it Поступила 27 октября 2014 г.}
\maketitle {\small
\begin{quote}
\noindent{\bf Аннотация. } В пространстве голоморфных функций в выпуклой области, изучается проблема кратной интерполяции посредством сумм рядов экспонент, сходящихся равномерно на всех компактах в области. Дискретное множество узлов кратной интерполяции лежит на вещественной оси в области и имеет единственную конечную предельную точку. Получен критерий разрешимости этой проблемы в терминах распределения предельных  направлений  показателей экспонент в бесконечности. 
\medskip

\noindent{\bf Ключевые слова:}{ голоморфная функция, выпуклая область, кратная интерполяция, ряд экспонент, замкнутый идеал, замкнутый подмодуль, сильно сопряженное пространство, двойственность}
\medskip
\end{quote}
\begin{quote}
\noindent{\bf Abstract.} In the space of all holomorphic functions in a convex domain it is studied interpolation problem with multiplicities by means of sums of the series of exponentials, converging  uniformly on compact subsets of the domain. The discrete set of the interpolation nodes with multiplicities counted is located on the real axis in the domain and has only one finite limit point. It is obtained а criterion for solvability of the problem in the terms of distribution limit directions of exponents of exponentials at infinity.  
\medskip

\end{quote} }
\section{Формулировка задачи и предварительные сведения}

Пусть $D$ -- выпуклая область в $\mathbb C.$ Обозначим через $H(D)$ пространство голоморфных функций в $D$ с топологией равномерной сходимости на компактных подмножествах из $D$. 
Рассмотрим произвольное бесконечное дискретное в $\mathbb C$ множество комплексных чисел $\Lambda=\{\lambda_{n}\}_{n\in\mathbb N}.$ 

Обозначим $$ \Sigma({\Lambda, D})=\{f\in H(D): f(z)=\sum_{n=1}^{\infty} c_{n} e^{\lambda_nz},\,z\in D\}.$$ Предполагается, что сходимость ряда экспонент абсолютная для каждой точки $z\in D,$ а тогда  (\cite{Leonu}) такой ряд сходится в топологии пространства $H(D).$ Для многомерной ситуации это показано, например, в работе \cite{Mcont}.

Предположим, что  $D\cap\mathbb R$ не пустое множество. Пусть задано бесконечное дискретное в области $ D$ множество вещественных узлов интерполяции, $\mathcal M =\{\mu_{k}\}_{k=1}^{\infty},$  $\mathcal M\subset D\cap\mathbb R.$ Кроме того, будем полагать, что каждому узлу $\mu_{ k}\in\mathcal M$ приписана кратность $ m_{k}\in\mathbb N.$
Если $f,g\in H(D),$ будем писать $f\cong g$ на $\mathcal M,$ если $f^{(j)}(\mu_{k})=g^{(j)}(\mu_{k})$ для всех $k\in \mathbb N$ и $j=0, 1, \dots,m_k-1.$
 \medskip

Рассмотрим в $H(D)$ следующую проблему интерполяции с вещественными узлами посредством сумм рядов экспонент:

{\it  \noindent Для произвольного множества узлов $\mathcal M\subset D\cap\mathbb R$ и для любой функции $g\in H(D)$ существует функция $f\in \Sigma({\Lambda, D}),$ такая, что  $f\cong g$ на $\mathcal M.$}

В силу классического результата об интерполяции голоморфными функциями (\cite{Her}, Следствие 1.5.4), эта задача может быть сформулирована в традиционных терминах: 

{\it \noindent Для любых интерполяционных данных  $b^{j}_{ k}\in \mathbb C, \, k\in \mathbb N,\, j=0, 1, \dots,m_k-1,$ существует функция $f\in \Sigma({\Lambda, D}),$ такая что $f^{(j)}(\mu_{ k})=b^{j}_{k},$} для всех $k$ и $j.$ \medskip
 
Обозначим через $\psi_\mathcal{ M}$ функцию из $H(D)$ с нулями во всех узлах $ \mu_{k} \in\mathcal{ M},$ с кратностями $ m_{k}, $ и только в них, и определим 
$$\bigl(\psi_\mathcal{ M}\bigr)=\{h\in H({D}): h=\psi_\mathcal{ M}\cdot
r, \, r\in H({D})\}\leqno{(1)}$$ --- замкнутый идеал в $H({D}),$ порожденный функцией $\psi_\mathcal{M}.$ Легко видеть, что $\bigl(\psi_\mathcal{ M}\bigr)=I_\mathcal{ M}=\{h\in H({D}): h\cong 0\ \mbox{на}\ \mathcal{ M} \}$ 
  
Разрешимость  проблемы интерполяции в $H(D)$ суммами рядов экспонент с показателями $ \Lambda $ и множеством узлов $ \mathcal{M} $ равносильна существованию следующего представления: 
$$H(D)= \Sigma({\Lambda, D})+ \bigl(\psi_{\mathcal{M}}\bigr).\leqno{(2)}$$

Единственности интерполяции в условиях рассматриваемой задачи быть не может, то есть $\Sigma({\Lambda, D})\cap\bigl(\psi_{\mathcal{M}}\bigr)\neq\{0\}.$ Это доказано в работе \cite{MP}  для пространства целых функций, но приведенное там доказательство с очевидными изменениями переносится на рассматриваемый  случай. \medskip
 
{\noindent \it Если имеется представление (2) и $\Sigma({\Lambda, D})\subset X \subset H(D),$ то справедливо и представление $H(D)= X+ \bigl(\psi_{\mathcal{M}}\bigr).$}

 В работе \cite{Na}, в случае, когда $D=\mathbb C,$ а $X$ есть ядро некоторого оператора свертки в пространстве целых функций $H(\mathbb C),$ найдены достаточные условия для интерполяции функциями из ядра оператора свертки в терминах расположения нулей $\Lambda$ характеристической функции этого оператора. Множество $\mathcal{M}$ в \cite{Na} имеет две предельных точки $ \pm\infty.$ В работе \cite{MP} удалось найти другие методы доказательства и, для всех возможных случаев расположения предельных точек $\mathcal{M},$ были получены критерии разрешимости проблемы кратной интерполяции в $H(\mathbb C)$ посредством сумм рядов экспонент из $ \Sigma({\Lambda, \mathbb C})\subset X.$  В случае, когда множество узлов имеет две предельные точки $ \pm \infty,$ критерий в \cite{MP} формулируется в тех же терминах, как и в работе \cite{Na}.
 
В данной статье метод доказательства достаточности из \cite{MP} распространяется  на случай пространства голоморфных функций в выпуклой области. Получен критерий интерполяции в случае, когда $\mathcal{M}$ имеет единственную предельную точку на границе $\partial\,D$ области $D.$ Критерий состоит в том, что распределение предельных направлений показателей $\Lambda$ в бесконечности связывается с геометрической структурой части границы выпуклой области $D$, которая содержит эту предельную точку.

Доказательство достаточности условия сводится к рассмотрению интерполяции рядами экспонент в пространстве функций, голоморфных в некоторой полуплоскости. Кроме того, оказалось, что, в рассмотриваемой здесь ситуации, при доказательстве необходимости этого условия,  нужно использовать идеи совсем другой природы, по сравнению с пространством $H(\mathbb C).$ Дело в  том, что ряды экспонент,  абсолютно сходящиеся на некотором множестве, обладают свойством распространения сходимости (см. например,  \cite{Mcont}). Следует отметить, что аналитическое продолжение для элементов общих инвариантных подпространств, допускающих спектральный синтез, изучалось в работе  \cite{Kri3}. 
\section{Двойственная формулировка проблемы интерполяции. Схема сведения к интерполяции в ядре оператора свертки}
 В дальнейшем, как и в работе \cite{MP}, используется схема доказательства, описанная в работе \cite{NapP}, основанная на двойственности с использованием преобразования Лапласа $\mathcal L$  функционалов.  При доказательстве достаточности условий интерполяции, предлагается   рассматривать естественные двойственные утверждения, отдельно для каждого из возможных вариантов расположения предельных точек множества $\mathcal M.$ 

 Кратко опишем эту схему, так как в работе  \cite{MP} она описана достаточно подробно для пространства $H(\mathbb C).$ В случае $H(D)$ укажем на некоторые изменения в приведенных там рассуждениях. 
 
Обозначим через $P_D$--- пространство всех целых функций экспоненциального типа с традиционной топологией индуктивного предела, которая обеспечивает топологический изоморфизм между сильным сопряженным пространством $H^{*}(D)$ и пространством $P_D$, реализующийся с помощью  преобразования Лапласа $\mathcal L$  функционалов $F\in H^{*}(D).$ Точнее, линейное непрерывное взаимнооднозначное преобразование Лапласа $\mathcal L$ функционалов $F\in H^{*}(D)$ определяется следующим образом: $\mathcal L: F\longmapsto \mathcal L{F}(z)=\left\langle F_\lambda,e^{\lambda z}\right\rangle,\, \mathcal L{F}\in P_D.$
  
Топология в ${(LN}^{*})$-пространстве $P_D$ не описывается в терминах сходимости последовательностей, однако секвенциально замкнутые подпространства являются замкнутыми (\cite{Seb}). Точное определение сходимости последовательностей в этой топологии будет приведено в  доказательстве достаточности условий леммы 4.\medskip
  
Определим раздельно непрерывную билинейную форму $\left[\cdot, \cdot\right]: H(D)\times P_D\longmapsto \mathbb C,$ согласно формуле $\left[ \psi, \varphi\right]=\left<{\mathcal L}^{-1}\varphi,\,\psi\right>,\,\psi\in H(D),\,\varphi \in P_D.$
 С помощью отображения $\varphi \longmapsto \left[\cdot,\phi \right]= \left<{\mathcal L}^{-1}\varphi,\cdot \right>,$ где $ {\mathcal L}^{-1}\varphi  \in H^{*}(D),$  задается изоморфизм между $P_D$ и сильным сопряженным пространством $H^{*}(D).$ Согласно введенной двойственности, любая функция из пространства ${P_D}$ взаимнооднозначно соответствует некоторому линейному непрерывному функционалу из  $H^*(D).$
 
 Хорошо известно, что каждая функция $ G \in P_{D},\,G\not\equiv0,$ имеющая минимальный тип при порядке один, порождает в пространстве $H(D)$ оператор свертки $M_G: H(D)\longmapsto H(D),$ который в рассматриваемой двойственности можно определить как
$$ M_G [\psi](z)=\left[S_z\bigl(\psi(\lambda)\bigr),  G_\lambda\right]=\left<({\mathcal L}^{-1}G)_\lambda, \psi (z+\lambda)\right>,$$
где $ S_z $ --- оператор сдвига: $ S_z\bigl(\psi(\lambda)\bigr)=\psi(\lambda+z).$ 

Известно, что $M_G $ линейный, непрерывный и сюръективный оператор. Cопряженный оператор к оператору свертки $M_{G}$ это оператор  $ A_G $ умножения на характеристическую функцию $G,$ корректно определенный на функциях ${\omega \in P_D}$ следующим образом: $\omega \longmapsto G\cdot\omega$ (Подробности в  \cite {NapC}, \cite {Leon}).

Обозначим $\Ker M_{G}=\{f\in H(D): M_G[f]=0\}$
 --- ядро оператора свертки $ M_{G},$ которое является замкнутым подпространством в $H(D),$ инвариантным относительно оператора дифференцирования.
  
Подпространство $\Ker M_{G}$ допускает спектральный синтез (\cite {Leon}, \cite {Kra}), то есть совпадает с замыканием в топологии пространства $H(D)$ линейной оболочки множества всех полиномиально - экспоненциальных мономов $z^{\nu}e^{\lambda_nz},$ содержащихся в нем.\medskip
 
Подпространство рядов экспонент $\Sigma({\Lambda, D}),$ вообще говоря, не замкнутое в $H(D)$.
В связи с этим, в доказательстве достаточности условий интерполяции, для каждого из возможных вариантов расположения множества узлов $\mathcal M,$  выделяется подпоследовательность ${\widetilde\Lambda}$ из ${\Lambda},$ таким образом, чтобы она являлась нулевым множеством некоторой целой функции  $G\in P_D$ минимального типа, причем $\Ker M_{G}= \Sigma({\widetilde\Lambda, D}).$ Затем будет получено представление (2) с заменой ${\Lambda}$ на ${\widetilde\Lambda},$ но тогда оно будет справедливо и для ${\Lambda}.$\medskip

 После того, как выделена подпоследовательность ${\widetilde\Lambda},$ достаточно доказать следующие два утверждения.
 
 \noindent{\it  $\it (I)$ Подпространство $\Ker M_{G}+ \bigl(\psi_\mathcal{ M}\bigr)$ --- всюду плотное в пространстве $H(D);$}  

\noindent {\it $\it (II)$ Подпространство $\Ker M_{G}+ \bigl(\psi_\mathcal{ M}\bigr)$ --- замкнутое в пространстве $H(D).$ } \medskip

Замкнутый идеал $\bigl(\psi_\mathcal{ M}\bigr)$ определен выше в (1).
В дальнейшем в этом параграфе для упрощения обозначений $\psi=\psi_\mathcal{ M}.$ 

Если $X_1$--- подпространство в топологическом векторном пространстве $X$, через $ X_1^0$ обозначим его поляру (или аннулятор), то есть множество функционалов из $X^*,$ которые обращаются в нуль на $X_1.$ 

Утверждение {\it  $\it (I)$} равносильно тому, что $\bigl(\Ker M_{G}+ (\psi)\bigr)^0=\bigl(\Ker M_{G}\bigr)^0\cap\bigl((\psi)\bigr)^0=\{0\}.$ Из Леммы 2 работы \cite{Mer1} следует, что Утверждение {\it  $\it (II)$} равносильно тому, что подпространство $\bigl(\Ker M_{G}\bigr)^0+ \bigl((\psi)\bigr)^0$  --- замкнутое в $P_D.$  \medskip

Пространство $P_D$ -- модуль над кольцом многочленов. С учетом двойственности, поляра $\bigl(\Ker M_{G}\bigr)^0$ совпадает с подмодулем, определяемым как
$$\bigl(G\bigr)_{P_D}=\{h\in P_D:  h= G\cdot r;\,r\in P_D\}.\leqno{(3)}$$

В доказательстве достаточности в лемме 4 будет доказано, что $\bigl(G\bigr)_{P_D}=\bigl(G\bigr)\cap {P_D},$ где $\bigl(G\bigr)$ -- замкнутый идеал в $H(\mathbb C),$ порожденный функцией $G.$ В частности отсюда следует, что подмодуль $\bigl(G\bigr)_{P_D}$ -- замкнутый. \medskip

Как известно, $(M^{*})$-пространство $H(D)$ -- рефлексивное (\cite{Seb},  \cite{NapC}), то есть его сильное второе сопряженное пространство $H^{**}(D) $ канонически изоморфно пространству $H(D).$
Поэтому отображение $\psi \longmapsto \left[\psi, \cdot\right],$ с учетом этого канонического изоморфизма, определяет изоморфизм между $(M^{*})$-пространством $H(D)$ и  сильным сопряженным  $P_D^{*}.$ Любая функция из $H(D)$ взаимнооднозначно соответствует некоторому линейному непрерывному функционалу из сильно сопряженного пространства ${P_D}.$ 

Более точно, это отображение понимается следующим образом: канонический изоморфизм $ H(D)$ и $H^{**}(D) $ имеет вид $\psi \longmapsto\Theta\psi= F_\psi,\,F_\psi\in P_D^*,\, $  $\left<F_\psi, \varphi\right>= \left[ \psi, \varphi\right]=\left<{\mathcal L}^{-1}\varphi,\,\psi\right>. $
Здесь $\psi\in H(D),\,\varphi\in P_D.$ 

Каждая функция $\psi\in H(D),\,\psi\not\equiv0, $ порождает в пространстве целых функций экспоненциального типа $P_D$ оператор свертки $\widetilde M_\psi: P_D\longmapsto P_D,$ $\widetilde M_\psi[\varphi](z)=\left[(\Theta\psi)_\lambda, S_z\bigl(\varphi(\lambda)\bigr)\right], $ где $ S_z $ --- оператор сдвига, $ S_z\bigl(\varphi(\lambda)\bigr)=\varphi(\lambda+z),\,\lambda \in\mathbb C.$

Далее получаем, что $\widetilde M_\psi[\varphi](z)=\bigl<{(\mathcal L}^{-1}S_z\varphi)_\lambda ,\,\psi(\lambda)\bigr>=\bigl<{e^{z\lambda} (\mathcal L}^{-1}\varphi)_\lambda,\,\psi(\lambda)\bigr>= \bigl<{ (\mathcal L}^{-1}\varphi)_\lambda,\,e^{z\lambda}\psi(\lambda)\bigr>, \,\varphi\in P_D.$ 

Отметим, что, используя известную формулу для обратного преобразования Бореля (\cite{Leon1}), отсюда можно получить явное интегральное представление этого оператора (см., например, \cite{Na}, \cite{MP} ).

Известно, что $\widetilde M_\psi$  линейный, непрерывный и сюръективный оператор. Оператор $\widetilde M_\psi$ является сопряженным к оператору $\widetilde A_{\psi}$ умножения на функцию $\psi$ в пространстве $H(D),$ действующему на функциях $g\in H(D)$ следующим образом: $g \longmapsto \psi\cdot g.$ Оператор $\widetilde A_{\psi}$ -- линейный и  непрерывный, а его образ совпадает с замкнутым идеалом $(\psi).$  Обозначим $\Ker  \widetilde M_{\psi}=\{f\in P_D: \widetilde M_\psi[f]=0\}.$
 
С учетом двойственности, поляра $\bigl( (\psi)\bigr)^0$ совпадает с $\Ker \widetilde M_{\psi}.$\medskip

В начале этого параграфа была описана схема того, как доказательство существования представления (2) сводится к утверждениям {\it $\it (I)$} и {\it $\it (II)$}, а затем было доказано следующее.

{\noindent\bf Предложение 1.} {\it Утверждения {\it $\it (I)$} и {\it $\it (II)$}, в $(M^{*})$-пространстве $H(D),$ равносильны двум  двойственным утверждениям в $({LN}^{*})$-пространстве ${P_D},$ соответственно: 

\noindent $\it (I^{*})$ Справедливо равенство ${\bigl(G\bigr)_{P_D}}\cap \Ker \widetilde M_{\psi}=\{0\}.$

\noindent $(\it II^{*})$  Подпространство ${\bigl(G\bigr)_{P_D}}+\Ker \widetilde M_{\psi}$ --- замкнутое в пространстве ${P_D}.$}
 \section{ Вспомогательные результаты}
 
 Справедливо следующее простое, но важное, утверждение.
 
{\noindent \bf Предложение 2:}  {\it Пусть $D_1 $ -- некоторая область, причем $ D\subset D_1, $ и области имеют общие части границы, на которых лежат все предельные точки множества $\mathcal{ M}.$ Если найдены некоторые условия на $\Lambda,$ при выполнении которых имеет место представление (2) c множеством узлов $\mathcal{ M}$ для пространства $ H(D_1),$  тогда такое представление  имеется и для $ H(D),$} с тем же самым множеством узлов.

{\bf Доказательство.} Для любой функции $ g\in H(D)$ существует $g_1\in H(D_1),\,$ $g\cong g_1\mbox{на} \ \mathcal{ M}.$ Тогда $ g=g_1+(g-g_1)$ в области $D.$
По условию существует $f_1\in \Sigma(\Lambda, D_1) \subset\Sigma(\Lambda, D),$ такая, что $f_1\cong g_1\vert_\mathcal{ M}.$ В области $D$ получили представление $ g=f_1+(g_1-f_1)+(g-g_1). $ Функции в скобках лежат в $ H(D)$ и равны нулю на $\mathcal{ M}$ с учетом кратностей. Доказательство закончено.\medskip
  
В дальнейшем нам потребуются некоторые свойства полиномов из экспонент с вещественными показателями. Такие полиномы изучены в монографии  \cite{Bel}.
  
Рассмотрим произвольный полином из экспонент  вида 
$$p(z)=\sum_{k=0}^s a_k(z)e^{\omega_k z},\,\omega_0<\omega_1<\cdots<\omega_s,\leqno{(4)}$$ где $ a_k(z) $ --- некоторые многочлены, и пусть $ a_0\cdot a_s\not\equiv 0.$ 

Из Теоремы 12.9 монографии \cite{Bel} легко получить, что  справедливо следующее.
\begin{lemma}
Существует такое $c_1>0,$ что во внешности круга $\{z\in\mathbb C:\left|z\right|\geqslant c_1\}$ выполнено: существуют положительные постоянные $ c_2,\,c_3$ и два вещественных числа $ m_0,\, m_s, $ причем $m_0>m_s$ или $m_0=m_s=0,$ такие, что
$$\left|p(z)\right|\geqslant c_2e^{\omega_0 \Re z},\leqno{(5)}$$
для всех $ z $ в области $ U_0 =\{z\in\mathbb C:\Re(z+m_0\ln z)<-c_3\},$  и
$$\left|p(z)\right|\geqslant c_2e^{\omega_s \Re z},\leqno{(6)}$$
для всех $ z $ в области $ U_s =\{z\in\mathbb C:\Re(z+m_s\ln z)>c_3\}. $
\end{lemma}
\medskip

Для любого фиксированного $ c\in\mathbb R,$ рассмотрим кривую $ \Re(z+m\ln z)=c,\,m\neq0.$ Она симметрична относительно вещественной оси. Для $ m>0$ эта кривая лежит в некоторой полуплоскости $ \Re z < A,\,A>0,$ а для $ m<0$ она лежит в некоторой полуплоскости $ \Re z > -A,\,A>0,$. Если точка $ z=x+iy $ лежит на кривой, то $ \left|\dfrac{y}{x}\right|\to\infty, $ $ \arg z \to \dfrac{\pi}{2},$ $ |z|=|y|(1+o(1)),$ при $  \left|z\right|\to\infty.$ Рассматриваемая кривая асимптотически приближается к показательной кривой $ x+m\ln|y|=c.$\medskip

Зафиксируем $ \beta \in\Big[0,\dfrac{\pi}{2}\Big)$ и для $\alpha\in \Big[0,\dfrac{\pi}{2}-\beta\Big),$ обозначим   $A_\alpha(\beta)=\{z\in \mathbb C: |\arg z-\beta| \leqslant \alpha \}.$

\begin{lemma}
Пусть $ \omega_s<0.$ Для произвольного полинома из экспонент $ p $ вида (4), существует такое $ r=r(p)>0,$ что для всех $z,\,|z|>r,$ в угле $A_\alpha(\beta)$ справедлива следующая оценка 
$$|p(z)|\geqslant c_3e^{\omega_s\cos(\beta+\alpha)|z|}.\leqno{(7)}$$
 
\end{lemma}
{\bf Доказательство.}  Легко видеть, что все точки $ z $ из угла $A_\alpha(\beta),$ лежащие вне некоторого круга,  лежат в области  $U_s.$ Поэтому, из оценки (6) полинома из экспонент $ p $ в области $U_s$  для $ |z|>c_1 ,$ вытекает оценка вне некоторого круга  $|z|>r$ в угле $ A_\alpha(\beta).$ Неравенство (7) вытекает из (6): если $ z=|z|e^{i\varphi},$ то в этом угле $0>\omega_s\Re e^{i\varphi}\geqslant\omega_s\cos(\beta+\alpha).$ Лемма 2 доказана.

Пусть некоторая  выпуклая область $ D $ содержит все показатели $ \omega_k,\,k=0, 1,\dots, s, $ полинома $ p $ из экспонент вида (4). Тогда легко показать, что $ p\in P_D.$ \medskip

Следующая лемма по существу доказана в \cite{MP} в несколько другой формулировке.
Рассмотрим произвольную бесконечную дискретную последовательность комплексных чисел $\mathcal{V}=\{v_{j}\},$ причем $\Re v_{j}>0.$  Предположим, что
$$\limsup_ {j\to\infty}\dfrac{\Re v_{j}}{\ln\left|v_{j}\right|}=\infty.\leqno{(8)}$$
Для дальнейшего важно заметить, что если последовательность $\mathcal{V}$ лежат в угле $A_\alpha(\beta),$ то условие (8) выполняется.

Обозначим $I_{\mathcal{V}}= \{f\in H(\mathbb C): f(v_{j})=0,\,j\in\mathbb N\} $ -- 
замкнутый идеал в $H(\mathbb C).$ (сравните с (3)). 

\begin{lemma}
 В описанной ситуации, если для $\mathcal{V}$ выполнено условие (8), то никакой многочлен  из экспонент $ p\not\equiv 0 $ вида (4) не может содержаться в идеале $I_{\mathcal{V}}.$
  \end{lemma}
 Дело в том, что условие (8) на нули идеала противоречит оценке (6).

  \section{Основные результаты}
  
  Для множества $ \Lambda $ введем множество $ P(\Lambda)$ предельных направлений в бесконечности как совокупность точек $ s\in\mathbb S,$ для которых найдется последовательность $ \{ \lambda_{n_k}\}_{k\in\mathbb N}, $ что 
$ \lim_{k\to\infty }\lambda_{n_k}/|\lambda_{n_k}|=s,\, \lim_{k\to\infty} |\lambda_{n_k}|=\infty.$
Множество $ P(\Lambda)$ замкнутое.\medskip 

Аналоги следующего понятия, под различными названиями, часто возникают в комплексном анализе, например, при изучении эффекта аналитического продолжения сумм рядов экспонент, их аналогов, а также элементов инвариантных подпространств (см., например, \cite{Mcont}, \cite{Krio}, \cite{Krio1}, \cite{Kri3}).  
Приведем некоторые  необходимые нам определения и результаты из работы \cite{Mcont}:

 Обозначим $ \mathbb S=\{s\in \mathbb C: |s|=1\}.$ Пусть $ S $ -- замкнутое подмножество $ \mathbb S.$
Пусть $D$ -- некоторая область в $\mathbb C.$ Обозначим $h(\varphi)= \sup_{\sigma\in D}\Re (e^{i\varphi}\sigma).$ Если $ k(\varphi): \mathbb C\mapsto (-\infty,+\infty]$ -- опорная функция (в смысле $ \mathbb R^2 $) выпуклой области $ D,$ то $h(\varphi)=k(-\varphi).$

Легко также видеть, что $ h(\varphi) $ это опорная функция (в смысле $ \mathbb R^2 $) области, комплексно сопряженной с $ D.$  Функция $ H(z)=\sup_{\sigma\in D}\Re (z\sigma)=h(\varphi)|z|,\, z=|z|e^{i\varphi}\in \mathbb C,$ -- позитивно однородная, полунепрерывная снизу, выпуклая функция. Из этого несложно получить, что функция $h(\varphi)$ -- полунепрерывная снизу на $\mathbb S.  $

Для $s\in  \mathbb S,\,s=e^{i\varphi},$ определим функцию  $d(s)=k(-\varphi),\,d:\mathbb S\mapsto (-\infty,+\infty].$ Видим, что для каждого $ s\in \mathbb S,\,$  $d(s)$ -- это наименьшая верхняя грань проекций области $D$ на направление $\overline s=e^{-i\varphi}.$ 
 
Например, для заданных $t\in\mathbb S$ и $c\in\mathbb R,$ обозначим $\Pi_c (\overline t)=\{z\in\mathbb C: \Re(tz)< c\}$ --  полуплоскость с направлением $ \overline{t}$ внешней нормали к границе, точка $ z=c\overline t$ лежит на ее границе. Для $ D= \Pi_c (\overline t)$ имеем, что $d(s)=+\infty,\, s\not=t,$ и $d(s)=c,\,s=t.$ \medskip

Множество 
$$
 D_S =\{z\in\mathbb C: \Re(sz)<d(s),\,s=e^{i\varphi}\in S,\}
$$ называетcя $ S$-выпуклой оболочкой области $ D.$

По определению, $S$-выпуклая оболочка  $ D_{S} $ любой области $ D$ это пересечение, по всем $s=e^{i\varphi}\in S,$  множеств $\Pi(\overline{s}, D)=\{z\in\mathbb C: \Re(sz)<d(s)\}.$ Если существует $ t\in S:$  $ d(t)=\infty,$ тогда $\Pi(\overline{t}, D)=\mathbb C.$ Если при этом существует хотя бы одно число $ s\in S, $ для которого $ d(s)<\infty,$ такие $t\in S$ в определении $ D_{S} $ можно не учитывать. 

Если  $ d(s)<\infty,$ множество $\Pi(\overline{s}, D)$ -- это опорная полуплоскость области $ D,$ то есть $ D\subset \Pi(\overline{s}, D) $ и $ \partial\,D\cap\partial\,\Pi(\overline{s}, D)\not =\emptyset.$ Легко видеть, что $\Pi(\overline{s}, D)=\Pi_0(\overline{s})+\overline {s}d(s).$ Здесь $ \Pi_0(\overline{s})=\{z\in\mathbb C: \Re(sz)<0\}).$ 

$D_S $ -- выпуклая область, более того, она $ S$-выпуклая (\cite{Mcont}, \cite{Krio}).
 Если $S=\mathbb S$, $S$-выпуклая оболочка множества  -- это обычная выпуклая оболочка.\medskip

  {\noindent \it {\bf Предложение А.}  Пусть $ D $ выпуклая область и $ S= P(\Lambda).$ Если ряд экспонент $\sum_{n=1}^{\infty}c_n e^{\lambda_nz}$ абсолютно сходится для всех $z\in D,$ то он абсолютно сходится и для $z\in D_{P(\Lambda)}.$ Его сумма - аналитическая функция в выпуклой области $ D_{P(\Lambda)}.$}\medskip

 Первое утверждение вытекает из предложений 16 и 8 работы \cite{Mcont}. В работе   \cite{Leonu} доказано, что ряд, абсолютно сходящий в выпуклой области $ D,$ сходится и в топологии пространства $ H(D) $ равномерной сходимости на компактах. \medskip

  {\bf Область  $ D $ -- полуплоскость.}   
Зафиксируем $\beta,\,|\arg\beta|<\dfrac{\pi}{2},$ и рассмотрим случай, когда $ D=\Pi_0(e^{-i\beta})$ -- "левая" полуплоскость. Обозначим $  s_\beta =e^{i\beta},$ тогда $ D=\Pi_0(\overline{s_\beta}). $ 
 
Пусть задано произвольное бесконечное дискретное множество вещественных узлов интерполяции $\mathcal M\subset \Pi_0(\overline{s_\beta})\cap \mathbb R^-.$  Каждая точка $ \mu_k\in \mathcal M$ имеет кратность $m_k, \, m_k\in \mathbb N.$ 
\medskip
 
\begin{lemma}
Пусть множество $ \mathcal {M}$ имеет единственную предельную точку $ z=0.$  В пространстве $ H\bigl(\Pi_0(\overline{s_\beta})\bigr)$ разрешима проблема кратной интерполяции рядами экспонент из $ \Sigma\bigl(\Lambda,\Pi_0(\overline{s_\beta})\bigr)$ с множеством узлов $\mathcal { M},$ тогда, и только тогда, когда $s_\beta \in P(\Lambda).$
\end{lemma} \medskip
Отметим, что направление $s_{\beta}$ -- комплексно сопряженное к направлению $\overline {s_\beta}=e^{-i\beta}$ внешней нормали  к границе $\partial\, \Pi_0(\overline{s_\beta}).$

{\bf Доказательство.} Из соображений симметрии в доказательстве можно считать, что $\beta\in \Big[0,\dfrac{\pi}{2}\Big).$ Условие леммы означает, что множества $\Lambda\cap A_\alpha(\beta)$ -- бесконечные, для всех достаточно малых $\alpha$. 

{\bf Необходимость.}
Пусть проблема кратной интерполяции рядами экспонент из $ \Sigma\bigl(\Lambda,\Pi_0(\overline{s_\beta})\bigr)$ с множеством узлов $ \mathcal { M} $ разрешима. Предположим, что $ s_\beta \not\in P(\Lambda).$ Тогда замкнутое множество $P(\Lambda)$ отделено от  направления  $s_\beta ,$ сопряженного к направлению $\overline{s_\beta} ,$ внешней нормали к границе $\partial\, \Pi_0(\overline{s_\beta}).$

Для любого числа $ s=e^{i\varphi}\in P(\Lambda)$ выполнено $ s\not=s_\beta.$ Поэтому, для области $ D= \Pi_0(\overline{s_\beta}),$ $d(s)=+\infty.$ Для любого $ s\in P(\Lambda)$ множество $\Pi(\overline{s}, D)=\mathbb C.$ 
По определению $ S$-выпуклой оболочки, $ D_{P(\Lambda)}= \mathbb C. $

Из предложения А получаем следующее. Если ряд экспонент абсолютно сходится в $\Pi_0(\overline{s_\beta})$ и $ s_\beta \not\in P(\Lambda)$, тогда этот ряд абсолютно сходится всюду в $\mathbb C.$ Тогда его сумма -- целая функция. 

Интерполяция целыми функциями с произвольными (например, неограниченными) данными  на множестве узлов $\mathcal { M},$ имеющем конечную предельную точку, невозможна. Противоречие.
\medskip

{\bf Достаточность.} Доказательство состоит из двух этапов.

1. Сначала сведем задачу к интерполяции в ядре некоторого оператора свертки.  
Если  утверждение леммы доказано для $\widetilde \Lambda\subset\Lambda,$ то оно будут доказано и для $ \Lambda.$  В дальнейшем,  мы перейдем к специальному  подпространству в $ \Sigma(\Lambda, D),$ замкнутому в $ H(D).$ Для этого заменим множество показателей на некоторую подпоследовательность из $\Lambda$.

Переходя к подпоследовательности, можно считать, что: 1) $ \Lambda\subset A_\alpha(\beta),$ для некоторого малого $ \alpha  ,$ 2) $ P(\Lambda)=\{s_\beta\},$ и 3) выполняется условие разделенности 
$$|\lambda_{n+1}|>2|\lambda_{n}|. \leqno{(9)}$$
 
Обозначим через $G$ целую функцию с простыми нулями $\lambda_{n},$ 

 $$G(z)=\prod_{n=1}^{\infty} \left( 1-\dfrac{z}{\lambda_{n}} \right).$$
 
Величина $\delta= \limsup_{n\to\infty}\dfrac{1}{|\lambda_{n}|}\ln\dfrac{1}{\bigl|G^\prime(\lambda_{n})\bigl|}$ 
 есть индекс Гельфонда-Леонтьева. 
 
Из условия (9) вытекает, что функция $ G $ имеет минимальный тип при порядке 1, и индекс конденсации $ \delta= 0.$  Это показано в работе \cite{MP}.

Из результатов монографии  \cite{Leon} (Теорема 4.2.2)  вытекает следующее утверждение. 

{\it  Пусть $ \delta=0.$ Тогда любая функция из замыкания в топологии $ H(D) $  линейной оболочки системы полиномиально-экспоненциальных мономов с множеством показателей, имеющим конечную верхнюю плотность с учетом кратностей, представляется в виде ряда экспонент.}\medskip

Подпространство $\Ker M_{G}$ допускает спектральный синтез. 
Тогда, с учетом Теоремы 4.2.3 из монографии \cite{Leon}, получаем следующее утверждение. 

\noindent {\bf Предложение Б.} {\it \noindent  Ядро $\Ker M_{G}$ состоит из всех функций $f(z),$ которые представляются рядами экспонент,
$$f(z)=\sum_{n=1}^\infty c_{n} e^{\lambda_{n}z},\,z\in\mathbb C,$$ \noindent сходящимися в топологии пространства $H(D),$ то есть $ \Ker M_{G}= \Sigma(\Lambda, D).$} 
\medskip

Следует отметить, что в многомерном случае, в более общей ситуации инвариантных подпространств, в работе \cite{Kri1} изучался фундаментальный принцип (в нашей ситуации это утверждение предложения Б).  Самая общая постановка этой задачи для комплексной плоскости рассмотрена в  \cite{Kri4}. В работе \cite{Kri2} подробно изучен случай рядов с вещественными показателями $ \Lambda.$

В этих работах введена новая характеристика $ S_ \Lambda,$ используя которую удалось получить критерии наличия фундаментального принципа для инвариантных подпространств в выпуклых областях. В силу этапа 1, $ \Lambda$ лежит в угле, а тогда, повторяя почти дословно доказательство из \cite{Kri2}, с. 100, получаем, что $ S_ \Lambda=0,$ и предложение Б можно получить и из результатов (\cite{Kri1}, \cite{Kri4},  \cite{Kri2}).\medskip  

2. На этом этапе переходим к доказательству двойственных утверждений. 

Обозначим через $\psi$ -- функцию из $ H(D) $ с нулевым множеством $ \mathcal {M}$ с учетом кратностей $ m_k.$

В силу Предложения 1 разрешимость интерполяционной проблемы вытекает из двух двойственных утверждений:

\noindent {\it $\it (I^{*})$} Справедливо равенство ${\bigl(G\bigr)_{P_D}}\cap \Ker \widetilde M_{\psi}=\{0\}.$

\noindent {\it$( II^{*})$}  Подпространство $\bigl(G\bigr)_{P_D}+\Ker \widetilde M_{\psi}$ --- замкнутое в пространстве ${P_D}.$

Подмодуль  $\bigl(G\bigr)_{P_D}$ определен выше в (3). \medskip

Важным моментом в доказательстве утверждений $(I^{*})$ и $(II^{*})$ является следующий  известный факт. 

{\it Подпространство $\Ker \widetilde M_{\psi}\subset P_D$ представляет собой линейную оболочку системы всех мономов вида $\{z^\nu e^{\mu_{k} z}\},\, k\in\mathbb N,\,\nu=0, 1, \cdots,m_{k}-1,$ то есть оно состоит только из полиномов из экспонент вида (4), где $ \omega_k=\mu_k $.}
Это несложно доказываемый фундаментальный принцип для $\Ker \widetilde M_{\psi}$ в пространстве $ P_D.$
\medskip

Двойственное утверждение {\it $(I^{*})$} следует из Леммы 3: покажем, что полином из экспонент $p\in\Ker \widetilde M_{\psi},\,p\not\equiv0,$ не может принадлежать  $\bigl(G\bigr)_{P_D}.$  Действительно, после этапа 1, считаем, что множество $ \Lambda$ лежит в $A_\alpha(\beta).$  Из этого следует, что для последовательности $ v_k=\lambda_k $ выполнено условие (8) из Леммы 3. Заметим, что $I_{ \Lambda}=\bigl(G\bigr)$ -- замкнутый идеал в $ H(\mathbb C). $ Как уже отмечалось выше в (3), $\bigl(G\bigr)_{P_D}=\bigl(G\bigr)\cap P_D.$ Утверждение $(I^{*})$ доказано.

Докажем последнее равенство. По определению, $\bigl(G\bigr)_{P_D}\subset I_{ \Lambda}\cap P_D.$ Далее,  согласно теореме (\cite{Leon1}) о делении на функцию минимального типа в пространстве $ P_D,$ верно и обратное включение. Кроме того, подмодуль в правой части последнего равенства --- замкнутый,  так как топология в $ P_\mathbb C $ сильнее топологии поточечной сходимости. Значит подмодуль $(G)_{P_D}$ замкнут в $ P_D.$ Эти факты были необходимы выше, при выводе двойственной формулировки проблемы интерполяции (Предложение 1).  \medskip

Получили, что имеется алгебраическая прямая сумма $(G)_{P_D}\oplus\Ker \widetilde M_{\psi}.$ Докажем замкнутость этого подпространства в $ P_D$ (это утверждение $(II^{*})$). 
Как известно (\cite{Seb}), в $({LN}^{*})$ --- пространстве $P_D$ замкнутость любого подпространства $X$ равносильна его секвенциальной замкнутости.
\medskip

{\it Сходимость последовательности $\{g_l\}_{l\in \mathbb N}$ в $(LN^{*})$-топологии пространства $P_D$ означает следующее: 

\noindent 1. Последовательность $\{g_l\}$ сходится к $ g $ в топологии пространства $ H(\mathbb C);$

\noindent 2. Существуют такие $A>0,\,j\in\mathbb N, $ что для всех $ l\in\mathbb N $ справедлива оценка $$|g_l(z)|\leqslant Ae^{H_j(z)},\, z\in\mathbb C. \leqno {(10)}$$ }

Здесь $\{K_j\}$ -- произвольное фиксированное счетное исчерпание области $ D$  выпуклыми компактами: $K_j \subset \int K_{j+1}$ и $ D= \bigcup_{j\in \mathbb N}K_j, $ $ H_j(z)=\sup_{\sigma\in K_j} \Re z\sigma. $ Если $ z=|z|e^{i\varphi},$ $h_j(\varphi)= H_j(z)/|z|$ -- опорная функция (в смысле $ \mathbb R^2 $) компакта, комплексно сопряженного с $ K_j.$ 
\medskip

Рассмотрим произвольную последовательность $\{g_l\}_{l\in\mathbb N}$ функций из 
$(G)_{P_D}\oplus\Ker \widetilde M_{\psi}$ и предположим, что она сходится в пространстве $P_D$ к функции $g \in {P_D}.$ Покажем, что предельная функция $g$ принадлежит $(G)_{P_D}\oplus\Ker \widetilde M_{\psi}.$ 

  Последовательность $\{g_l\}$ состоит из функций вида $ g_l= p_l+R_l, $ где функции  $ R_l \in  (G)_{P_D},$ то есть $ R_l\vert_{\Lambda}=0,$ а функции $p_l\in\Ker \widetilde M_{\psi}.$
  
Если в последовательности  $\{g_l\}$ содержится бесконечно много членов с $ R_l \equiv 0,$ то предельная функция $ g\in \Ker \widetilde M_{\psi}.$ Если в $\{g_l\}$  содержится бесконечно много членов с $ p_l\equiv 0, $ то $ g\in  (G)_{P_D}.$  Для всех таких последовательностей $\{g_l\}$ предельная функция $g\in(G)_{P_D}\oplus\Ker \widetilde M_{\psi}.$ 

Следовательно, далее можно считать, что последовательность $\{g_l\}$ такова, что  $R_l\not \equiv 0,\,p_l\not\equiv 0$ для всех $ l.$\medskip

Полагаем, что $\mu_k<\mu_{k+1}<0,\,\mu_k\to0,\,k\to\infty.$
 Так как $p_l \in \Ker  \widetilde M_{\psi},\,  p\not\equiv 0, $ это полином из экспонент вида $$p_l(z)=\sum_{\Fin_{\mathcal {M}}^{(l)}} a_k^{(l)}(z)e^{\mu_k z}.$$
 Здесь, для любого $ k\in\mathbb N,$ функции $a_{k}^l$  --- произвольные многочлены степеней не выше $m_{k}-1,$ соответственно.
 Для каждого $ l\in \mathbb N $ справа стоит сумма по некоторому конечному подмножеству  $Fin_{\mathcal M}^l\subset\mathcal{M}.$ Обозначим через $u_l$ номер максимального из $ \mu_{k} $ в этом представлении, то есть $a_{u_l}^l \not\equiv 0. $ \medskip

  Пусть последовательность $\{g_l\}$ такова, что множество  чисел  $\{u_l\}$ бесконечное. Покажем, что оно ограниченное.
Предположим, что множество $\{u_l\}$ является неограниченным. 

Выберем в качестве исчерпания полуплоскости $ \Pi_0(\overline{s_\beta})$ полукруги $ K_j=e^{-i\beta}\cdot  B_j^-,$  где $ B_j^-=(-1/j+\{|z|\leqslant j\})\cap \{\Re z\leqslant -1/j\}. $ Для каждого $ j $ обозначим $ t_j= \art j^2.$ Обозначим $\varepsilon_j= \dfrac{1}{2}(\dfrac{\pi}{2}- t_j),$ тогда несложно показать, что справедливы оценки:
$$-\dfrac{1}{j}|z| \leqslant H_j(z) \leqslant-A_j|z|\  \mbox{для} \  z\in A_{\varepsilon_j}({\beta}).\leqno {(11)}$$ где $A_j=\dfrac{\sqrt{1+j^4}}{j^2}\sin \varepsilon_j =\dfrac{1}{2j^2}(1+o(j)),\,j\to\infty.$
\medskip
  
Все многочлены из экспонент $ p_l $ имеют вид (1). Выберем $ k>j,$ такое, что $\varepsilon_k \in  [0,\dfrac{\pi}{2}-\beta).$ Так как $ a_{-q_l}^{(l)}\not\equiv 0,$  можно применить  оценку (7) для $ p_l $  из леммы 2, в которой $ \alpha=\varepsilon_k.$ Используя еще оценку (10), получаем следующую оценку для $ R_l = g_l - p_l,$  $p_l\not \equiv 0,\,R_l\not \equiv 0 ,$ :  
$$ |R_l(z)|\geqslant |p_l(z)|-|g_l(z)|\geqslant c_3e^{\mu_{u_{l}}\cos(\beta+\alpha)|z|}-Ae^{H_j(z)},$$ для всех   $z$ в области $ \{z\in A_{\varepsilon_k}(\beta),\,|z|>r\}.$ Здесь $r=r(l).$ 
Так так $ k>j,\,A_{\varepsilon_k}(\beta)\subset A_{\varepsilon_j}(\beta),$ и тогда из (11) следует, что
$$ |R_l(z)|\geqslant |p_l(z)|-|g_l(z)|\geqslant c_3e^{\mu_{u_{l}}\cos(\beta+\alpha)|z|}-Ae^{-A_j|z| },$$
вне некоторого круга $ |z|>r $ в угле $ A_{\varepsilon_k}(\beta).$

По предположению, множество $u_l$ неограниченное, поэтому в представлениях полиномов $p_l$ из экспонент существуют  $\mu_{u_l},$ сколь угодно близкие к 0.

Выберем $\mu_{u_{l_0}}>A_j/\cos(\beta+\alpha),$ тогда из последней оценки вытекает, что $ |R_{l_0}(z)|>0$ для всех $z$ вне некоторого круга $\{|z|>r_1(l_0)\}$ в угле $ A_{\varepsilon_k}(\beta).$

Получили противоречие: действительно, в силу этапа 1, $ P(\Lambda)=\{s_\beta\},$ поэтому, для любого $ k $ вне любого круга, в угле $ A_{\varepsilon_k}(\beta)$ лежит бесконечная последовательность точек из $\Lambda,$ а мы предполагали, что $ R_{l_0} \vert_{\Lambda}=0.$

Следует отметить, что для произвольной последовательности $ \{g_l\}$  компакт $K_j$ может быть сколько угодно большим, а тогда величина $ \varepsilon_k$ может быть сколь угодно малой. В этом смысл условия леммы. \medskip 

Итак, в представлениях полиномов $p_l$ из экспонент в произвольной сходящейся последовательности  $\{g_l\},$ $g_l=p_l+R_{l},$ множество чисел $u_{l}$  ограниченное.
Следовательно, последовательность $\{p_l\}$ принадлежит некоторому конечномерному подпространству $X\subset \Ker  \widetilde M_{\psi_1}.$ Утверждение $(I^{*})$ означает, что все элементы сходящейся последовательности $ g_l= p_l+R_l$ лежат в алгебраической прямой сумме $X\oplus (G)_{P_D}\subset \Ker  \widetilde M_{\psi_1}\oplus (G)_{P_D}.$ 

В любом топологическом векторном пространстве алгебраическая сумма конечномерного подпространства и замкнутого подпространства является замкнутым подпространством (\cite{Rud}, стр. 41). Итак, предельная функция $g$  последовательности  $ g_l= p_l+R_l$  принадлежит $\Ker  \widetilde M_{\psi_1}\oplus (G)_{P_D}.$ Утверждение $(II^{*})$ доказано.

Из доказанных утверждений $(I^{*})$ и $(II^{*})$ вытекает утверждение леммы 4.
\medskip

{\noindent \bf Замечание 1.} В доказательстве достаточности показано следующее. Пусть $ \Lambda$ -- произвольное множество показателей. Тогда одно лишь общее условие (8) (оно следует из условий леммы 4) -- достаточное для того, чтобы множество $ \Sigma({\Lambda, \Pi_{-\beta}})+I_\mathcal{M}$ было всюду плотным в топологии пространства $ H(\Pi_{-\beta}).$ 

{\noindent \bf Замечание 2.} После преобразования $z\to-z$ плоскости $ \mathbb C,$ получим формулировку, соответствующую случаю "правой" полуплоскости. Кроме того, для любого $ h\in \mathbb R, $ в рассматриваемой задаче допустимо преобразование $z\to z+h$ комплексной плоскости, после которого соответствующим образом нужно изменить формулировки. Действительно, в результате этого преобразования, множество рядов экспонент сохраняется, а множество узлов сдвигается.
\medskip

 {\bf Область $ D$ -- выпуклая, на ее границе лежит конечная предельная точка $ \mathcal M $.}
 Пусть $D$ -- выпуклая область в $\mathbb C.$
 
Обозначим $h(\varphi)= \sup_{\sigma\in D}\Re (e^{i\varphi}\sigma).$ Для каждого $ \varphi,\,$ число $ h(\varphi)$ -- это значение опорной функции $ k(-\varphi) $ (в смысле $ \mathbb R^2$) области $ D$ в направлении $ e^{-i\varphi }.$ Пусть $ s= e^{i\varphi },$  ранее была определена функция $ d(s)=k(-\varphi).$

Прямая $ l(\overline s)=\{z=x+iy: \Re(sz)=x\cos (-\varphi)+y\sin(-\varphi)=d\}$ называется опорной для области $D$ в направлении $\overline s=e^{ -i\varphi},$ если на границе $D $  существует точка, принадлежащая  $ l(\overline s),$ причем область $D$ лежит в опорной полуплоскости $\{z\in \mathbb C: \Re(sz)<d\}.$  Эту точку назовем точкой опоры для прямой $ l(\overline s).$ Легко видеть, что прямая $ l(\overline s)$ -- опорная, тогда, и только тогда, когда $d=d(s).$ В силу выпуклости, $ D\subset \Pi(\overline{s}, D).$ \medskip
 
Пусть $ 0\in \partial\,D.$ Обозначим через $T_D(0)\subset\mathbb S$ совокупность всех $ s\in\mathbb S,$ для которых точка $0$ на границе $ D$ является точкой опоры для $ l(\overline{s}).$ Ясно, что $T_D(0)=\{s\in\mathbb S:  d(s)=0\}.$  

Заметим, что, в условиях леммы 4, $ D=\Pi_0(\overline{s_\beta}),$ $T_D(0)=\{s_\beta\}.$

 \begin{Theorem 1} 
 Пусть $ D$ выпуклая область, причем  $ 0\in \partial\,D $ и $ D\cap\mathbb R \not = \emptyset.$  Предположим, что множество $\mathcal M\subset D\cap\mathbb R$ -- дискретное в $ D$ и имеет единственную предельную точку $z= 0.$ В пространстве $ H(D)$ разрешима проблема кратной интерполяции рядами экспонент из $ \Sigma(\Lambda, D)$ с множеством узлов $\mathcal { M},$ тогда, и только тогда, когда множество $  P(\Lambda)\cap T_D(0)\not=\emptyset.$ 
\end{Theorem 1}

{\noindent\bf Доказательство.} Случай $ D=\Pi_0(\overline{s_\beta}),\,|\beta|<\dfrac{\pi}{2},$ рассмотрен в лемме 4.

Без ограничения общности можно считать, что $D\cap\mathbb R^- \not = \emptyset,$ тогда $ \mathcal {M}\subset D \cap\mathbb R^-.$ В противном случае можно использовать преобразование $ z\to-z $ плоскости  $ \mathbb C$. 

Из того, что $ h(\varphi)\geqslant 0,$ и полунепрерывности снизу следует, что множество $T_D(0)$ -- замкнутое. Из выпуклости и однородности функции $H(z)=h(\varphi)|z| $ следует, что $T_D(0)$ связное множество.
 Легко также видеть, что $|\arg s|<\dfrac{\pi}{2} $ для всех $s\in T_D(0),$ так как $ \mathcal {M}\cap \mathbb R^-\not=\emptyset.$

Необходимость условия $  P(\Lambda)\cap T_D(0)\not=\emptyset$  следует из Предложения Б, а достаточность доказывается сведением к лемме 4.

{\noindent\bf Необходимость.}
Предположим, что проблема интерполяции разрешима, но условие теоремы не выполнено,  $  P(\Lambda)\cap T_D(0)=\emptyset$. Далее будет показано, что в этом случае точка $ z=0 $ лежит в $D_{P(\Lambda)}.$

Множества $ P(\Lambda)$ и $T_D(0)=\{s\in\mathbb S:  d( s)=0\}$ замкнутые. 

Для любого подмножества $ X\subset \mathbb S $ и числа $\delta>0$ обозначим$ X_\delta=\{s\in\mathbb S:\exists u\in  X,\ |s-u|\leqslant\delta\}.$

Найдется такое $ \delta >0,$ что $ P(\Lambda)\cap\bigl( (T_D(0)\bigr)_\delta=\emptyset,$ поэтому существует такое связное замкнутое множество $ S_1\in \mathbb S ,$ что $ P(\Lambda)\in \int S_1,\,S_1\cap T_D(0)=\emptyset.$ 

Тогда, из определения $ S$-выпуклой оболочки вытекает, что 
$$D_{S_1}\subset D_{P(\Lambda)}. \leqno (12)$$

Для всех $s\in S_1$ выполнено  $ d(s)>0,$ поэтому из полунепрерывности следует, что $\exists c,\, d(s)>c>0,\, s\in S_1.$ Обозначим $ B(c)=\{z\in \mathbb C: |z|=c\}.$ Для всех $ z\in B(c)$ и любого $s\in S_1,$ $ \Re(sz)<c< d(s).$ Отсюда следует, что точка $0\in \partial\,D$ лежит в $\Pi_c(\overline s)\subset \Pi(\overline s, D),$ для любого $ s\in S_1.$ Доказано, что $0\in\big(B(c)\big)_{S_1}\subset D_{S_1}.$ 

Из (12) получаем, что $0\in D_{P(\Lambda)}.$ 
Доказательство завершается следующим образом. 

В силу предложения А, любой ряд экспонент, который  сходится абсолютно в выпуклой области $ D,$ абсолютно сходится и в выпуклой области $ D_{P(\Lambda)}.$ Его сумма - аналитическая функция в $ D_{P(\Lambda)}.$

Точка $ z=0$ лежит в области $ D_{P(\Lambda)}$ и, по условию, она является предельной для множества узлов $\mathcal {M}.$ Интерполяция аналитическими функциями для произвольных (например, неограниченных) данных в множестве узлов, имеюшем предельную точку в области, невозможна. Противоречие. Доказательство необходимости условия теоремы закончено.

{\noindent \bf Достаточность.}
По условию, существует предельное направление $s_\beta=e^{i\beta} \in P(\Lambda),$ лежащее в $T_D(0),$ тогда $|\beta|<\dfrac{\pi}{2},$ так как множество $D\cap\mathbb R^-$ непустое.

Так как $ s_\beta \in T_D(0),$ точка $ z=0$ -- точка опоры прямой $ l(\overline{s_\beta})=\{z: \Re (zs_\beta)=0\}.$ Область $ D $ лежит в опорной полуплоскости $ \Pi_0(\overline{s_\beta})=\{\Re (zs_\beta)<0\},$ и $ 0\in \partial\,\Pi_0(\overline{s_\beta})= l(\overline{s_\beta}).$ Таким образом, множество $\mathcal { M}\subset \Pi_0(\overline{s_\beta})$ можно использовать в качестве множества узлов для интерполяции рядами экспонент в пространстве $ H\bigl(\Pi_0(\overline{s_\beta})\bigr)\subset H(D).$  При этих условиях, разрешимость проблемы интерполяции рядами экспонент в пространстве $ H\bigl(\Pi_0(\overline{s_\beta})\bigr) $  доказана в лемме 4. Разрешимость проблемы в $ H(D)$ следует из Предложения 2. Доказательство закончено.\medskip



\end{document}